\documentclass[a4paper,12pt]{amsart}

\usepackage[active]{srcltx}
\usepackage{amsfonts}
\usepackage{amssymb}
\usepackage{ifthen}
\usepackage{graphicx}
\newtheorem{thm}{Theorem}
\newtheorem{cor}[equation]{Corollary}
\newtheorem{lem}{Lemma}
\newtheorem{prop}[equation]{Proposition}

\newtheorem{conj}[equation]{Conjecture}
\newtheorem{rem}{Remark}
\theoremstyle{definition}
\newtheorem{defn}{Definition}%[section]

\newtheorem{prob}[equation]{Problem}
\newtheorem{ques}[equation]{Question}

\newcounter {own}
\def\theown {\thesection       .\arabic{own}}

\newenvironment{pf}[1][]{%
 \vskip 3mm
 \noindent
 \ifthenelse{\equal{#1}{}}%
  {{\slshape Proof. }}%
  {{\slshape #1.} }%
 }%
{\qed\bigskip}

\newcounter{alphabet}

\newenvironment{Thm}[1][]{\refstepcounter{alphabet}%
\bigskip%
\noindent%
{\bf Theorem \Alph{alphabet}}%
\ifthenelse{\equal{#1}{}}{}{ (#1)}%
{\bf .} \itshape}{\vskip 8pt}

\newenvironment{Lem}[1][]{\refstepcounter{alphabet}%
\bigskip%
\noindent%
{\bf Lemma \Alph{alphabet}}%
{\bf .} \itshape}{\vskip 8pt}

\newcommand{\Aut}{{\operatorname{Aut}}}

\def\be{\begin{equation}}
\def\ee{\end{equation}}

\newcommand{\bee}{\begin{enumerate}}
\newcommand{\eee}{\end{enumerate}}

\newcommand{\blem}{\begin{lem}}
\newcommand{\elem}{\end{lem}}
\newcommand{\bthm}{\begin{thm}}
\newcommand{\ethm}{\end{thm}}
\newcommand{\bcor}{\begin{cor}}
\newcommand{\ecor}{\end{cor}}
\newcommand{\beg}{\begin{examp}}
\newcommand{\eeg}{\end{examp}}
\newcommand{\begs}{\begin{examples}}
\newcommand{\eegs}{\end{examples}}
\newcommand{\bdefe}{\begin{defn}}
\newcommand{\edefe}{\end{defn}}
\newcommand{\bprob}{\begin{prob}}
\newcommand{\eprob}{\end{prob}}
\newcommand{\bques}{\begin{ques}}
\newcommand{\eques}{\end{ques}}
\newcommand{\bei}{\begin{itemize}}
\newcommand{\eei}{\end{itemize}}

\newcommand{\bcon}{\begin{conj}}
\newcommand{\econ}{\end{conj}}
\newcommand{\bcons}{\begin{conjs}}
\newcommand{\econs}{\end{conjs}}
\newcommand{\bprop}{\begin{prop}}
\newcommand{\eprop}{\end{prop}}
\newcommand{\br}{\begin{rem}}
\newcommand{\er}{\end{rem}}
\newcommand{\brs}{\begin{rems}}
\newcommand{\ers}{\end{rems}}
\newcommand{\bo}{\begin{obser}}
\newcommand{\eo}{\end{obser}}
\newcommand{\bos}{\begin{obsers}}
\newcommand{\eos}{\end{obsers}}
\newcommand{\bpf}{\begin{pf}}
\newcommand{\epf}{\end{pf}}
\newcommand{\ba}{\begin{array}}
\newcommand{\ea}{\end{array}}
\newcommand{\beq}{\begin{eqnarray}}
\newcommand{\beqq}{\begin{eqnarray*}}
\newcommand{\eeq}{\end{eqnarray}}
\newcommand{\eeqq}{\end{eqnarray*}}

\newcommand{\ds}{\displaystyle}

\begin{document}
%\begin{center}
%{ Communicated by  }
%\end{center}
\bibliographystyle{amsplain}
\title [] {Lipschitz continuity  of holomorphic mappings with respect to Bergman metric}

\author{Shaolin Chen}
\address{S. L. Chen, Department of Mathematics, Soochow University, Suzhou, Jiangsu
215006, People's Republic of China;   College of Mathematics and
Statistics, Hengyang Normal University, Hengyang, Hunan 421008,
People's Republic of China.} \email{mathechen@126.com}

\author{David Kalaj}
\address{D. Kalaj, Faculty of Natural Sciences and Mathematics,
University of Montenegro, Cetinjski put b. b. 81000 Podgorica,
Montenegro. } \email{davidk@t-com.me}

%\author{X. Wang${}^{~\mathbf{*}}$}
%\address{X. Wang, Department of Mathematics,
%Hunan Normal University, Changsha, Hunan 410081, People's Republic
%of China.} \email{xtwang@hunnu.edu.cn}

\subjclass[2000]{Primary: 32A10; Secondary: 32A17, 32A35}
\keywords{ Holomorphic mapping, Lipschitz continuity, Bergman metric.}
%\\
%${}^{\mathbf{*}}$ Corresponding author}
\date{\today  %November 4, 10;
File: Lip.con.tex}
\begin{abstract}
In this paper, we establish the sharp estimate of the Lipschitz
continuity with respect to the Bergman metric. The obtained results
are the improvement and  generalization  of the corresponding
results of Ghatage, Yan and Zheng (Proc. Amer. Math. Soc., 129:
2039-2044, 2000).
\end{abstract}

%\thanks{The research was partly supported by
%NSF of China (No. 11071063). The work was carried out while the
%first author was visiting IIT Madras, under ``RTFDCS Fellowship."
%This author thanks Centre for International Co-operation in Science
%(Formerly  Centre for Cooperation in Science \& Technology
%among Developing Societies (CCSTDS)) for its support and cooperation.
%}

\maketitle \pagestyle{myheadings} \markboth{S. L. Chen and D.
Kalaj}{Lipschitz continuity  of holomorphic mappings with respect to
Bergman metric}

\section{Introduction and main results}\label{csw-sec1}

Let $\mathbb{C}^{n}$ denote the Euclidean space of complex dimension
$n$. Throughout this paper, we write a point $z\in\mathbb{C}^{n}$ as
a column vector in $n\times1$ matrix form
$z=(z_{1},\ldots,z_{n})^{T},$ where the symbol $T$ stands for the
transpose of vectors or matrices. For $z\in \mathbb{C}^{n}$, the
conjugate of $z$, denoted by $\overline{z}$, is defined by
$\overline{z}=(\overline{z}_{1},\ldots, \overline{z}_{n} )^{T}. $
For $z\in\mathbb{C}^{n}$ and
$w=(w_{1},\ldots,w_{n})^{T}\in\mathbb{C}^{n} $, we write $\langle
z,w\rangle := \sum_{k=1}^nz_k\overline{w}_k$ and $ |z|:={\langle
z,z\rangle}^{1/2}=(|z_{1}|^{2}+\cdots+|z_{n}|^{2})^{ 1/2}. $ For
$a=(a_{1},\ldots,a_{n})^{T}\in \mathbb{C}^n$, we set
$\mathbb{B}^n(a, r)=\{z\in \mathbb{C}^{n}:\, |z-a|<r\}. $ Also, we
use $\mathbb{B}^n$ to denote the unit ball $\mathbb{B}^n(0, 1)$ and
let $\mathbb{D}=\mathbb{B}^1$.

The class of all holomorphic functions from $\mathbb{B}^{n}$ into
$\mathbb{C}^{n}$ is denoted by $H(\mathbb{B}^{n},\mathbb{C}^{n})$.
Let $\Aut(\mathbb{B}^{n})$ be the automorphism group consisting of
all biholomorphic self mappings of the unit ball $\mathbb{B}^{n}$.
 %We recall the following facts from \cite{R}:
%\begin{enumerate}
%\item[\textbf{(a)}] For $a\in\mathbb{B}^{n}$, let
%$$\phi_{a}(z)=\frac{a-P_{a}z-(1-|a|^{2})^{1/2} Q_{a}z}{1-\langle z,a\rangle},$$
%where $$P_{a}z=\frac{a\langle z, a\rangle}{\langle a, a\rangle}
%~\mbox{ and }~ Q_{a}z=z-P_{a}z.$$ Then
%$\phi_{a}\in\Aut(\mathbb{B}^{n})$.

%\item[\textbf{(b)}] For $z\in\mathbb{B}^{n}$ and
%$\phi\in\Aut(\mathbb{B}^{n})$, \be\label{eq1.04}
%|\phi'(z)\theta|\geq\frac{1-|\phi(z)|^{2}}{(1-|z|^{2})^{1/2}} \ee
%and \be\label{eq1.05}
%|\det\phi'(z)|=\left(\frac{1-|\phi(z)|^{2}}{1-|z|^{2}}\right)^{(n+1)/2},
%\ee where $\theta\in\partial\mathbb{B}^{n}.$
%\end{enumerate}

For $z\in\mathbb{B}^{n}$, let
$$B(z)=\frac{(1-|z|^{2})I+A(z)}{(1-|z|^{2})^{2}}$$ be
{\it the Bergman matrix}, where $I$ is the $n\times n$ identity
matrix and

$$A(z)=\left(\begin{array}{cccc}
\ds z_{1}\overline{z}_{1}\; \cdots\;
 z_{1}\overline{z}_{n}\\[4mm]
 %5\ds \frac{\partial v_{1}}{\partial x_{1}}\; \frac{\partial v_{1}}{\partial y_{1}}\;
%\frac{\partial v_{1}}{\partial x_{2}}\; \frac{\partial
%v_{1}}{\partial y_{2}}\;\cdots\;
 %\frac{\partial v_{1}}{\partial x_{n}}\; \frac{\partial v_{1}}{\partial
 %y_{n}}\\[2mm]
\vdots \\[2mm]
 \ds z_{n}\overline{z}_{1}\; \cdots\;
 z_{n}\overline{z}_{n}
\end{array}\right).
$$
For a smooth curve $\gamma:~[0,1]\rightarrow\mathbb{B}^{n}$, let

$$\ell(\gamma)=\int_{0}^{1}\langle B(\gamma(t))\gamma'(t),\gamma'(t)\rangle^{1/2}dt.$$
%A metric
%$\beta:\mathbb{B}^{n}\times\mathbb{B}^{n}\rightarrow[0,\infty)$ can
%be defined as follows.
 For any two points $z$ and $w$ in
$\mathbb{B}^{n}$, let $\beta(z,w)$ be the infimum of the set
consisting of all $\ell(\gamma)$, where $\gamma$ is a piecewise
smooth curve in $\mathbb{B}^{n}$ from $z$ to $w$. Then we call
$\beta$ {\it the Bergman metric} in $\mathbb{B}^{n}$ (cf. \cite{Z}).

As in \cite{LX}, the \textit{prenorm}
$\|f\|_{\mathcal{P}(n,\alpha)}$ of $f\in H(\mathbb{B}^{n},
\mathbb{C}^{n})$ is given by
$$\|f\|_{\mathcal{P}(n,\alpha)}=\sup_{z\in\mathbb{B}^{n}}
D_{f}^{n,\alpha}(z),
$$ where $\alpha>0$ and $$D_{f}^{n,\alpha}(z)=(1-|z|^{2})^{\frac{\alpha(n+1)}{2n}}|\det
f'(z)|^{\frac{1}{n}}.$$ Let $\mathcal{B}_{\mathcal{P}(n,\alpha)}$ be
the class of all holomorphic mappings $f\in H(\mathbb{B}^{n},
\mathbb{C}^{n})$ satisfying $\|f\|_{\mathcal{P}(n,\alpha)}<\infty$.
In particular,
 $\mathcal{B}_{\mathcal{P}(1,\alpha)}$ is the classical family of
$\alpha$-Bloch functions  (cf. \cite{Z1,Z}).

For $z,w\in\mathbb{D}$, the {\it pseudo-hyperbolic distance} is
defined as
$$\rho(z,w)=\left|\frac{z-w}{1-\overline{w}z}\right|.$$ In
particular, if $n=1$, then, for $z,w\in\mathbb{D}$,
$\tanh\beta(z,w)=\rho(z,w).$ In \cite{GYZ},  Ghatage, Yan and  Zheng
showed that $D_{f}^{1,1}(z)$ is Lipschitz continuous with respect to
the pseudo-hyperbolic metric, which is given in Theorem~A below. For more details
on this topic, see \cite{CG,CPW-1,MM,X}.

\begin{Thm}{\rm (\cite[Theorem 1]{GYZ})}\label{Thm-1}
Let $f\in\mathcal{B}_{\mathcal{P}(1,1)}$. Then, for all
$z,w\in\mathbb{D}$,

$$\big|D_{f}^{1,1}(z)-D_{f}^{1,1}(w)\big|\leq3.31\|f\|_{\mathcal{P}(1,1)}\rho(z,w).$$

\end{Thm}

Let $\varphi$ be a holomorphic mapping of $\mathbb{B}^{n}$ into
$\mathbb{B}^{n}$. For all $f\in\mathcal{B}_{\mathcal{P}(n,\alpha)}$,
let $C_{\varphi}:~ f\mapsto f\circ\varphi$ be a {\it
pre-composition} operator. As an application of Theorem~A,
they proved

\begin{Thm}{\rm (\cite[Theorem 2]{GYZ})}\label{Thm-2}
Let   $\varphi$ be a holomorphic mapping of $\mathbb{D}$ into
$\mathbb{D}$.   If for some constants $r\in(0,1/4)$, and
$\varepsilon>0$, for each $w\in\mathbb{D}$, there is a point
$z_{w}\in\mathbb{D}$ such that
$$\rho(\varphi(z_{w}),w)<r,~\mbox{and}~\frac{1-|z_{w}|^{2}}{1-|\varphi(z_{w})|^{2}}|\varphi'(z_{w})|>\varepsilon,$$
 then
$C_{\varphi}:~\mathcal{B}_{\mathcal{P}(1,1)}\rightarrow\mathcal{B}_{\mathcal{P}(1,1)}$
is bounded below.
\end{Thm}

In this paper, by using a different  method,   we generalize
Theorems~A and ~B to several dimensional case and
obtain the sharp estimate of the Lipschitz constant with respect
to the Bergman metric. 

\begin{thm}\label{thm-1} Let $f\in \mathcal{B}_{\mathcal{P}(n,1)}$. Then, for
$z_{1},z_{2}\in\mathbb{B}^{n}$,
%$$\left|(1-|z_{1}|^{2})^{\frac{(n+1)}{2n}}|\det f'(z_{1})|^{\frac{1}{n}}-
%(1-|z_{2}|^{2})^{\frac{(n+1)}{2n}}|\detf'(z_{2})|^{\frac{1}{n}}\right|\leq C\|f\|_{SN(1)}\tanh\beta(z_{1},z_{2}).$$
\be\label{eq-c1}\big|D^{n,1}_{f}(z_{2})-D^{n,1}_{f}(z_{1})\big|\leq
M(n)
\|f\|_{\mathcal{P}(n,1)}\big[\tanh\beta(z_{1},z_{2})\big]^{\frac{1}{n}},\ee
where
$M(n)=(2+n)^{\frac{1}{2n}}\left(\frac{n+2}{n+1}\right)^{\frac{n+1}{2n}}.$
Moreover, the constant $M(n)$ in {\rm (\ref{eq-c1})} cannot be
replaced by a smaller number.
\end{thm}

%We remark that if $n=1$ in Theorem \ref{thm-1}, then, for
%$z_{1},z_{2}\in\mathbb{B}^{n}$,  we have
%$$\big|D^{1,1}_{f}(z_{2})-D^{1,1}_{f}(z_{1})\big|\leq\frac{3\sqrt{3}}{2}
%\|f\|_{\mathcal{P}(1,1)}\rho(z_{1},z_{2}).$$

The following result is an application of Theorem \ref{thm-1}.

\begin{thm}\label{cor-1}
Let $f\in \mathcal{B}_{\mathcal{P}(1,1)}$. Then, for
$z\in\mathbb{D}$,

\be\label{eq-1.19}(1-|z|^{2})\left(\Big|\frac{\partial}{\partial
z}D^{1,1}_{f}(z)\Big|+\Big|\frac{\partial}{\partial
\overline{z}}D^{1,1}_{f}(z)\Big|\right)
\leq\frac{3\sqrt{3}}{2}\|f\|_{\mathcal{P}(1,1)}\ee and
\be\label{eq-1.20}|f''(z)|\leq\frac{\left(2|z|+\frac{3\sqrt{3}}{2}\right)\|f\|_{\mathcal{P}(1,1)}}{(1-|z|^{2})^{2}}.\ee

Moreover, the extreme functions $f(z)=\pm3\sqrt{3}z^{2}/4$ show that
the estimates {\rm (\ref{eq-1.19})} and {\rm (\ref{eq-1.20})} are
sharp.

%$$|f^{(n)}(z)|\leq\frac{\left(2+\frac{3\sqrt{3}}{2}\right)^{n-1}\|f\|_{\mathcal{P}(1,1)}}{1-|z|^{2}}.$$
\end{thm}

%\begin{rem} From Theorem {\rm \ref{cor-1}}, we see that if $f\in \mathcal{B}_{\mathcal{P}(1,1)}$, then, for any fixed
%$n$, $$|f^{(n)}(z)|=O\left(\frac{1}{1-|z|^{2}}\right)$$ as
%$|z|\rightarrow1^{-}$. This is a surprising phenomenon.\end{rem}

Applying Theorem \ref{thm-1}, we get the following result which is
an improvement of Theorem~B.

\begin{thm}\label{thm-2}
Let $\varphi$ be a holomorphic mapping of $\mathbb{B}^{n}$ into
$\mathbb{B}^{n}$. Suppose that there is constants
$0<r<\frac{1}{M(n)}\left(\frac{n+2}{1+n}\right)^{\frac{1}{n}}$ and
$\varepsilon>0$ such that, for  each $w\in\mathbb{B}^{n}$, there is
a point $z_{w}\in\mathbb{B}^{n}$ satisfying
$\tanh\beta(\varphi(z_{w}),w)<r^{n}$ and
$|\tau_{\varphi}(z_{w})|>\varepsilon$, where
$$\tau_{\varphi}(z_{w})=\left(\frac{1-|z_{w}|^{2}}{1-|\varphi(z_{w})|^{2}}\right)^{\frac{n+1}{2n}}
\left|\det \varphi'(z_{w})\right|^{\frac{1}{n}}$$ and $M(n)$ is
defined as in Theorem {\rm \ref{thm-1}}.
 Then, for all
$f\in\mathcal{B}_{\mathcal{P}(n,1)}$, there is a  constant
$k(n,r,\varepsilon)>0$ depended only on $r$, $\varepsilon$ and $n$
such that
$$\|C_{\varphi}(f)\|_{\mathcal{P}(n,1)}\geq k(n,r,\varepsilon)\|f\|_{\mathcal{P}(n,1)}.$$
\end{thm}

The proofs of Theorems \ref{thm-1}-\ref{thm-2} will be presented in
Section \ref{csw-sec2}.

\section{Proofs of the main results}\label{csw-sec2}

We begin the section by recalling the following results which play
an important role in the proofs of Theorem \ref{thm-1}.

\begin{Lem}{\rm (\cite[Lemma 1.1]{CPW-1})}\label{Lem1}
For $x\in[0,1]$, let
$$\varphi(x)=x(1-x^{2})^{\frac{\alpha(n+1)}{2}}\sqrt{\alpha(1+n)+1}\left [\frac{\alpha(n+1)+1}
{\alpha(n+1)}\right ]^{\frac{\alpha(n+1)}{2}}
$$
and
$$a_0(\alpha)=\frac{1}{\sqrt{\alpha(1+n)+1}}.
$$
Then $\varphi$ is increasing in $ [0,a_0(\alpha) ]$, decreasing in $
[a_0(\alpha) ,1]$ and $\varphi (a_0(\alpha))=1.$
\end{Lem}

\begin{Thm}{\rm (\cite[Theorem 1.2]{CPW-1})}\label{Thm1}
Suppose that $f\in H(\mathbb{B}^{n}, \mathbb{C}^{n})$ such that
$\|f\|_{\mathcal{P}(n,\alpha)}=1$ and $\det f'(0)=\lambda\in(0,1]$.
Then, for all $z$ with $|z|\leq
\frac{a_0(\alpha)+m_{\alpha}(\lambda)}{1+a_0(\alpha)m_{\alpha}(\lambda)}$,
we have \be\label{eq1.07} |\det f'(z)|\geq {\rm Re\,}\big(\det
f'(z)\big)\geq
\frac{\lambda(m_{\alpha}(\lambda)-|z|)}{m_{\alpha}(\lambda)(1-m_{\alpha}(\lambda)|z|)^{\alpha(n+1)+1}},
\ee where $m_{\alpha}(\lambda)$ is the unique real root of the
equation $\varphi(x)=\lambda$ in the interval $[0,a_0(\alpha)]$ and,
$\varphi$ and $a_0(\alpha)$ are defined as in Lemma~{\rm C}. Moreover, for all $z$ with
$|z|\leq\frac{a_0(\alpha)-m_{\alpha}(\lambda)}{1-a_0(\alpha)m_{\alpha}(\lambda)}$,
we have \be\label{eq1.08} |\det
f'(z)|\leq\frac{\lambda(m_{\alpha}(\lambda)+|z|)}{m_{\alpha}(\lambda)(1+m_{\alpha}(\lambda)|z|)^{\alpha(n+1)+1}}.
\ee Moreover, the estimates of $(\ref{eq1.07})$  and
$(\ref{eq1.08})$ are sharp.
\end{Thm}

\subsection*{Proof of Theorem \ref{thm-1}} Without loss of generality, we assume that
$\|f\|_{\mathcal{P}(n,1)}=1$ and $B^{n,1}_{f}(z_{2})\leq
B^{n,1}_{f}(z_{1}).$ Let $\phi\in\Aut(\mathbb{B}^{n})$ such that
$\phi(0)=z_{1}$ and $w=\phi^{-1}(z_{2})$. For $z\in\mathbb{B}^{n}$,
set $g=f(\phi(z))$. By \cite[Proposition 1.21]{Z}, we have

$$\tanh\beta(z_{1},z_{2})=\tanh\beta\big(\phi^{-1}(z_{1}),\phi^{-1}(z_{2})\big)=\tanh\beta(0,w)=|w|.$$
Since
$$|\det\phi'(z)|=\left(\frac{1-|\phi(z)|^{2}}{1-|z|^{2}}\right)^{\frac{n+1}{2}},$$
we see that \be\label{eq-1.10}|\det g'(0)|^{\frac{1}{n}}=|\det
f'(\phi(0))|^{\frac{1}{n}}|\det \phi'(0)|^{\frac{1}{n}}=|\det
f'(z_{1})|^{\frac{1}{n}}(1-|z_{1}|^{2})^{\frac{n+1}{2n}}\ee and

%\begin{eqnarray*}
\beq\label{eq-1.11}(1-|w|^{2})^{\frac{n+1}{2n}}|\det
g'(w)|^{\frac{1}{n}}&=&(1-|w|^{2})^{\frac{n+1}{2n}}|\det
f'(\phi(w))|^{\frac{1}{n}}|\det \phi'(w)|^{\frac{1}{n}}\\ \nonumber
&=&(1-|w|^{2})^{\frac{n+1}{2n}}|\det
f'(z_{2})|^{\frac{1}{n}}\left(\frac{1-|\phi(w)|^{2}}{1-|w|^{2}}\right)^{\frac{n+1}{2n}}\\
\nonumber &=&|\det
f'(z_{2})|^{\frac{1}{n}}(1-|z_{2}|^{2})^{\frac{n+1}{2n}}.\eeq
%\end{eqnarray*}

$\mathbf{Case~1.}$ If $|\det g'(0)|=0$, then it is obvious.

$\mathbf{Case~2.}$ Let $\det g'(0)=\lambda e^{i\theta}$, where
$\lambda>0$ and $\theta\in[0,2\pi]$. Applying Theorem~D
(\ref{eq1.07}) to $e^{-i\frac{\theta}{n}}g(z)$, for
$|z|\leq\frac{a_{0}(1)+m_{1}(\lambda)}{1+a_{0}(1)m_{1}(\lambda)}$,
we have

\be\label{eq-1.12}\mbox{Re}\left(e^{-i\theta}\det g'(z)\right)\geq
\frac{\lambda(m_{1}(\lambda)-|z|)}{m_{1}(\lambda)(1-m_{1}(\lambda)|z|)^{n+2}},\ee
where $m_{1}(\lambda)$ and $a_0(1)$ are defined as in Theorem~
D.

$\mathbf{Subcase~2.1.}$  $|w|\leq m_{1}(\lambda)$. It is easy to
know that
\be\label{eq-cd-1}m_{1}(\lambda)\leq\frac{a_{0}(1)+m_{1}(\lambda)}{1+a_{0}(1)m_{1}(\lambda)}.\ee
On the other hand, by calculations, we have
$$(1-|w|^{2})^{\frac{n+1}{2n}}\geq(1-|w|m_{1}(\lambda))^{\frac{n+1}{2n}}\geq(1-|w|m_{1}(\lambda))^{\frac{n+2}{n}}$$ and
$$m_{1}^{\frac{1}{n}}(\lambda)-|w|^{\frac{1}{n}}\leq\big(m_{1}(\lambda)-|w|\big)^{\frac{1}{n}},$$
which gives \be\label{eq-1.13} m_{1}^{\frac{1}{n}}(\lambda)-
\frac{(1-|w|^{2})^{\frac{n+1}{2n}}(m_{1}(\lambda)-|w|)^{\frac{1}{n}}}
{(1-m_{1}(\lambda)|w|)^{\frac{n+2}{n}}}\leq|w|^{\frac{1}{n}}. \ee

By Lemma~C  and Theorem ~D, we have
\be\label{eq-1.14}m_{1}(\lambda)(1-m_{1}^{2}(\lambda))\sqrt{2+n}\left(\frac{n+2}{n+1}\right)^{\frac{n+1}{2}}=\lambda,\ee
which, together with (\ref{eq-1.10}), (\ref{eq-1.11}),
(\ref{eq-1.12}), (\ref{eq-cd-1}) and (\ref{eq-1.13}), implies that

\begin{eqnarray*}
D^{n,1}_{f}(z_{1})-D^{n,1}_{f}(z_{2})&=&|\det
g'(0)|^{\frac{1}{n}}-(1-|w|^{2})^{\frac{n+1}{2n}}|\det
g'(w)|^{\frac{1}{n}}\\&=&\lambda^{\frac{1}{n}}-(1-|w|^{2})^{\frac{n+1}{2n}}|\det
g'(w)|^{\frac{1}{n}}\\
&\leq&\lambda^{\frac{1}{n}}-\frac{(1-|w|^{2})^{\frac{n+1}{2n}}\lambda^{\frac{1}{n}}(m_{1}(\lambda)-|w|)^{\frac{1}{n}}}
{m_{1}^{\frac{1}{n}}(\lambda)(1-m_{1}(\lambda)|w|)^{\frac{n+2}{n}}}\\
&=&\left(\frac{\lambda}{m_{1}(\lambda)}\right)^{\frac{1}{n}}\left[m_{1}^{\frac{1}{n}}(\lambda)-
\frac{(1-|w|^{2})^{\frac{n+1}{2n}}(m_{1}(\lambda)-|w|)^{\frac{1}{n}}}
{(1-m_{1}(\lambda)|w|)^{\frac{n+2}{n}}}\right]\\
&\leq&\left(\frac{\lambda}{m_{1}(\lambda)}\right)^{\frac{1}{n}}|w|^{\frac{1}{n}}\\
&=&M(n)\big(1-m_{1}^{2}(\lambda)\big)^{\frac{1}{n}}|w|^{\frac{1}{n}}\\
&\leq&M(n)|w|^{\frac{1}{n}}.
\end{eqnarray*}

$\mathbf{Subcase~2.2.}$  $|w|> m_{1}(\lambda)$. Then, by
(\ref{eq-1.10}), (\ref{eq-1.11}) and (\ref{eq-1.14}),

\begin{eqnarray*}
D^{n,1}_{f}(z_{1})-D^{n,1}_{f}(z_{2})&=&|\det
g'(0)|^{\frac{1}{n}}-(1-|w|^{2})^{\frac{n+1}{2n}}|\det
g'(w)|^{\frac{1}{n}}\\&\leq&\lambda^{\frac{1}{n}}\\
&=&m_{1}^{\frac{1}{n}}(\lambda)\big(1-m_{1}^{2}(\lambda)\big)^{\frac{n+1}{2n}}M(n)\\
&<&M(n)|w|^{\frac{1}{n}}.
\end{eqnarray*}

Now we prove the sharpness part. For $z\in\mathbb{B}^{n}$, let
%$$f_{\lambda}(z)=\left(\begin{array}{cccc} \displaystyle
%\int_{0}^{z_{1}}\frac{\lambda(m_{1}(\lambda)-\xi)}{m_{1}(\lambda)(1-m_{1}(\lambda)\xi)^{n+2}}\,d\xi   \\
%z_2\\ \vdots \\ z_{n}\end{array}\right).$$

$$f_{\lambda}(z)=\left(\begin{array}{cccc}
\int_{0}^{z_{1}}\frac{\lambda(m_{1}(\lambda)-\xi)}{m_{1}(\lambda)(1-m_{1}(\lambda)\xi)^{n+2}}\,d\xi   \\
z_{2}\\
\vdots \\
 z_{n}
\end{array}\right).
$$
%$$f_{\lambda}(z)=\left(\int_{0}^{z_{1}}\frac{\lambda(m_{1}(\lambda)-\xi)}{m_{1}(\lambda)(1-m_{1}(\lambda)\xi)^{n+2}}\,d\xi, z_{2}, \ldots,z_{n}\right)'.$$
For any $\epsilon\in\left(0,M(n)\right]$, let
$$m_{1}(\lambda)=\min\left\{\bigg[1-\Big(1-\frac{\epsilon}{M(n)}\Big)^{\frac{2n}{n+1}}\bigg]^{\frac{1}{2}},
M(n)\right\}.$$ Then, for $w'=(0,0,\ldots,0)$ and
$w''=(m_{1}(\lambda),0,\ldots,0)$, we have

\begin{eqnarray*}
\left|D^{n,1}_{f_{\lambda}}(w')-D^{n,1}_{f_{\lambda}}(w'')\right|&=&\lambda^{\frac{1}{n}}\\
&=&M(n)m_{1}^{\frac{1}{n}}(\lambda)\left(1-m^{2}_{1}(\lambda)\right)^{\frac{1}{n}}\\
&>&\|f\|_{\mathcal{P}(n,1)}\big[\tanh\beta(w',w'')\big]^{\frac{1}{n}}\big(M(n)-\epsilon\big),
\end{eqnarray*}
which shows that the constant $M(n)$ is sharp. The proof of this
theorem is complete. \hfill $\Box$

\subsection*{Proof of Theorem \ref{cor-1}}
For $z=x+iy\in\mathbb{D}$, let $w=z+re^{i\theta}$. Then, by Theorem
\ref{thm-1}, we have

%\begin{eqnarray*}
\beq\label{eq-1.16} \nonumber
\Lambda_{f}(z)&=&\max_{\theta\in[0,2\pi]}\left[\lim_{r\rightarrow0^{+}}
\bigg(\frac{\big|D^{1,1}_{f}(z)-D^{1,1}_{f}(w)\big|}{r}\cdot\frac{r}{\rho(z,w)}\bigg)\right]\\
\nonumber
&=&(1-|z|^{2})\max_{\theta\in[0,2\pi]}\left|\frac{\partial}{\partial
x}D^{1,1}_{f}(z)\cos\theta+\frac{\partial}{\partial
y}D^{1,1}_{f}(z)\sin\theta\right|\\ \nonumber
&=&\frac{(1-|z|^{2})}{2}\left(\bigg|\frac{\partial}{\partial
x}D^{1,1}_{f}(z)+i\frac{\partial}{\partial
y}D^{1,1}_{f}(z)\bigg|+\bigg|\frac{\partial}{\partial
x}D^{1,1}_{f}(z)-i\frac{\partial}{\partial
y}D^{1,1}_{f}(z)\bigg|\right)\\ \nonumber
&=&(1-|z|^{2})\left(\Big|\frac{\partial}{\partial
z}D^{1,1}_{f}(z)\Big|+\Big|\frac{\partial}{\partial
\overline{z}}D^{1,1}_{f}(z)\Big|\right)\\
&\leq&\frac{3\sqrt{3}}{2}\|f\|_{\mathcal{P}(1,1)}, \eeq
%\end{eqnarray*}
where
$$\Lambda_{f}(z)=\max_{\theta\in[0,2\pi]}\left(\lim_{r\rightarrow0^{+}}\frac{\big|D^{1,1}_{f}(z)-D^{1,1}_{f}(w)\big|}{\rho(z,w)}\right).$$
On the other hand,
\begin{eqnarray*}
\Big|\frac{\partial}{\partial
z}D^{1,1}_{f}(z)\Big|+\Big|\frac{\partial}{\partial
\overline{z}}D^{1,1}_{f}(z)\Big|&=&\left|-\overline{z}|f'(z)|+\frac{f''(z)\overline{f'(z)}}{2|f'(z)|}(1-|z|^{2})\right|\\
&&+\left|-z|f'(z)|+\frac{\overline{f''(z)}f'(z)}{2|f'(z)|}(1-|z|^{2})\right|\\
&\geq&|f''(z)|(1-|z|^{2})-2|z||f'(z)|,
\end{eqnarray*}
which, together with (\ref{eq-1.16}), gives that
\begin{eqnarray*}
%\beq\label{eq-1.17}
|f''(z)|(1-|z|^{2})&\leq&2|z||f'(z)|+\Big|\frac{\partial}{\partial
z}D^{1,1}_{f}(z)\Big|+\Big|\frac{\partial}{\partial
\overline{z}}D^{1,1}_{f}(z)\Big|\\
&\leq&\frac{2|z|\|f\|_{\mathcal{P}(1,1)}}{1-|z|^{2}}+\frac{\frac{3\sqrt{3}}{2}\|f\|_{\mathcal{P}(1,1)}}{1-|z|^{2}}\\
&=&\frac{\left(2|z|+\frac{3\sqrt{3}}{2}\right)\|f\|_{\mathcal{P}(1,1)}}{(1-|z|^{2})}.
\end{eqnarray*} The proof of this theorem is complete. \hfill $\Box$
%\end{eqnarray*}
%Conclusions can be drawn by a similar calculation procedure

%By (\ref{eq-1.17}), we see that
%$$|f''(z)|(1-|z|^{2})\leq\left(2+\frac{3\sqrt{3}}{2}\right)\|f\|_{\mathcal{P}(1,1)},$$
%which, together with the similar calculation procedure of
%(\ref{eq-1.17}), yields that

%$$|f^{(n)}(z)|\leq\frac{\left(2+\frac{3\sqrt{3}}{2}\right)^{n-1}\|f\|_{\mathcal{P}(1,1)}}{1-|z|^{2}}.$$

\subsection*{Proof of Theorem \ref{thm-2}} Without loss of generality, we assume that
$\|f\|_{\mathcal{P}(n,1)}=1$. For $z\in\mathbb{B}^{n}$, it follows
from  Theorem \ref{thm-1} that there is a point $w\in\mathbb{B}^{n}$
such that
$$D_{f}^{n,1}(w)>1-\sigma$$ and
$$\big|D_{f}^{n,1}(w)-D_{f}^{n,1}(z) \big|\leq \left[M(n)\bigg(\frac{n+1}{n+2}\bigg)^{\frac{1}{n}}+
\sigma\right]\big[\tanh\beta(w,z)\big]^{\frac{1}{n}},$$ where
$$\sigma=\frac{1-rM(n)\Big(\frac{n+1}{n+2}\Big)^{\frac{1}{n}}}{2(1+r)}~
\mbox{and}~r<\frac{1}{M(n)}\left(\frac{n+2}{1+n}\right)^{\frac{1}{n}}.$$
By the assumption,  there is a point $z_{w}$ such that
$$\big[\tanh\beta(\varphi(z_{w}),w)\big]^{\frac{1}{n}}<r<\frac{1}{M(n)}\left(\frac{n+2}{1+n}\right)^{\frac{1}{n}}$$
and $|\tau_{\varphi}(z_{w})|>\varepsilon$, which imply that

%\begin{eqnarray*}
\beq\label{eq-1.15}\nonumber
D^{n,1}_{f_{\lambda}}(\varphi(z_{w}))&\geq&
D^{n,1}_{f_{\lambda}}(w)-\left[M(n)\bigg(\frac{n+1}{n+2}\bigg)^{\frac{1}{n}}+
\sigma\right]\big[\tanh\beta(\varphi(z_{w}),w)\big]^{\frac{1}{n}}\\
\nonumber
&\geq&1-\sigma-\left[M(n)\bigg(\frac{n+1}{n+2}\bigg)^{\frac{1}{n}}+
\sigma\right]r\\ \nonumber
&=&1-rM(n)\bigg(\frac{n+1}{n+2}\bigg)^{\frac{1}{n}}-(1+r)\sigma\\
%\nonumber
&=&\frac{1-rM(n)\Big(\frac{n+1}{n+2}\Big)^{\frac{1}{n}}}{2}>0. \eeq
%\end{eqnarray*}
By (\ref{eq-1.15}), we conclude that

$$\|C_{\varphi}(f)\|_{\mathcal{P}(n,1)}\geq D^{n,1}_{f_{\lambda}}(\varphi(z_{w}))
|\tau_{\varphi}(z_{w})|>k(n,r,\varepsilon),$$ where
$$k(n,r,\varepsilon)=\frac{\left[1-rM(n)\Big(\frac{n+1}{n+2}\Big)^{\frac{1}{n}}\right]\varepsilon}{2}.$$
The proof of this theorem is compete. \hfill $\Box$\\

{\bf Acknowledgement:} This research of the first author was partly
supported by the National Natural Science Foundation of China (No.
11401184 and No. 11571216),  the Construct Program of the Key
Discipline in Hunan Province,  the Science and Technology Plan of
Hunan Province (2016TP1020),  the Fifty-ninth Batch of Post Doctoral
Foundation of China (No. 2016M590492) and the Post Doctoral
Foundation of Jiangsu Province (No. 1601182C).

%1601182C

%{\bf Acknowledgement:} The authors thank the referee for bringing the work of \cite{B,E,Ra1,Ra2} to their attention.

\end{document}